\documentclass[12pt]{amsart}
\usepackage{bbding}
\usepackage[margin=2.0cm]{geometry}

\usepackage{amsfonts}
\usepackage{mathrsfs}
\usepackage{cite}
\usepackage{tikz}

\newcommand{\R}{{\mathbb R}}
\newcommand{\C}{{\mathbb C}}

\newtheorem{theorem}{Theorem}[section]
\newtheorem{lemma}[theorem]{Lemma}
\newtheorem{prop}[theorem]{Proposition}
\newtheorem{coro}[theorem]{Corollary}
\newtheorem{claim}[theorem]{Claim}

\newtheorem{example}[theorem]{Example}

\title{The planar $L_p$-Minkowski problem for $0< p<1$}
\author{K\'aroly J. B\"or\"oczky}
\address{Alfr\'ed R\'enyi Institute of Mathematics, Hungarian Academy
  of Sciences, Reltanoda u. 13-15, H-1053 Budapest, Hungary, and
Department of Mathematics, Central European University, Nador u 9, H-1051, Budapest, Hungary, boroczky.karoly.j@renyi.mta.hu}
\author{Hai T. Trinh}
\address{Department of Mathematics, Central European University, Nador u 9, H-1051, Budapest, Hungary, haitrinh1210@gmail.com}
\thanks{
Research is supported in parts by NKFIH grants 109789, 121649 and 116451}

\begin{document}
\maketitle

\begin{abstract}
Necessary and sufficient conditions for the
existence of solutions to the asymmetric $L_p$ Minkowski problem in $\R^2$ are
established for $0 < p < 1$.
\end{abstract}

\section{Introduction}

For the notions of the Brunn-Minkowski theory in $\R^n$ used in this paper, see Schneider \cite{SCH}. We write $\mathcal{H}^m$, $m\leq n$, for the $m$-dimensional Hausdorff measure normalized in a way such that it coincides with the Lebesgue measure on $\R^m$. We call a compact convex set $K$ with non-empty interior in $\R^n$ a convex body. For $x\in \partial K$, we choose an exterior unit normal $\nu_K(x)$ to $\partial K$ at $x$, which is unique for
$\mathcal{H}^{n-1}$ almost all $x\in\partial K$. 
The surface area measure $S_K$ on $S^{n-1}$ is defined for a Borel set $\omega\subset S^{n-1}$ by
$$
S_{K}(\omega)=\int_{x\in\nu_{K}^{-1}(\omega)}d\mathcal{H}^{n-1}(x).
$$
The classical Minkowski existence theorem, due to Minkowski in the case of polytopes or discrete measures and to Alexandrov for the general case, states that a Borel measure $\mu$ on $S^{n-1}$ is the surface area measure of a convex body if and only if the measure of any open hemisphere is positive, and
$$
\int_{S^{n-1}}ud\mu(u)=0.
$$
The solution is unique up to translation. If the measure $\mu$ has a density function $f$ with respect
to $\mathcal{H}^{n-1}$ on $S^{n-1}$, then even the regularity of the solution is well understood, see
Lewy \cite{LE}, Nirenberg \cite{NIR}, Cheng and Yau \cite{CY}, Pogorelov \cite{POG}, and Caffarelli \cite{LC}.

Lutwak  \cite{LUT} initiated the study of the so called $L_{p}$ surface area measure for any $p\in\R$. For a convex compact set $K$ in $\R^n$, let $h_K$ be its support function, and hence
$$
h_K(u)=\max\{\langle x,u\rangle:\, x\in K\} \mbox{ \ \ for $u\in\R^n$} 
$$
where $\langle\cdot,\cdot\rangle$ stands for the Euclidean scalar product.
Let ${\mathcal K}_0^n$ denote family of convex bodies in $\R^n$ containing the origin $o$.  For $p\leq 1$ and 
$K\in {\mathcal K}_0^n$, the  $L_{p}$-surface area measure is defined by
$$
dS_{K,p}=h_K^{1-p}\,d S_K.
$$
In particular, if  $\omega\subset S^{n-1}$ Borel, then
\begin{equation}
\label{SKpdef}
S_{K,p}(\omega)=\int_{x\in\nu_{K}^{-1}(\omega)}\langle x,\nu_{K}(x)\rangle^{1-p}d\mathcal{H}^{n-1}(x).
\end{equation}
Here the case $p=1$ corresponds to the surface area measure $S_K$, and $p=0$ to the so called cone volume measure.
If $p>1$, then the same formula $dS_{K,p}=h_K^{1-p}\,d S_K$ defines the $L_{p}$-surface area measure, only one needs to assume that  
 either $o\in{\rm int}\,K$, or $o\in\partial K$ and $\int_{S^{n-1}}h_K^{1-p}\,d S_K<\infty$.

The $L_{p}$ surface area measure has been intensively investigated in the recent decades, see, for example, \cite{ADA, BG, CG, GM, H1, HP1, HP2, HENK, LU1, LU2, LR, LYZ1, LYZ2, LYZ3, LYZ6, LZ, NAO, NR, PAO, PW}. In \cite{LUT}, Lutwak posed the associated $L_{p}$ Minkowski problem for $p\geq 1$ which extends the classical Minkowski problem.
If $p>1$ and $p\neq n$, then the $L_{p}$ Minkowski problem is solved by 
Chou, Wang \cite{CW},  Guan, Lin \cite{GL} and  Hug, Lutwak, Yang, Zhang \cite{HLYZ2}.
In addition, the $L_p$ Minkowski problem for $p<1$ was publicized by a series of talks by Lutwak in the 1990's. 
 The $L_{p}$ Minkowski problem is the classical Minkowski problem when $p=1$, while the $L_{p}$ Minkowski problem is the so called logarithmic Minkowski problem when $p=0$, see, for example,
\cite{BoH16, BLYZ, BLYZ2, BLYZ3, LU1, LU2, LR, NAO, NR, PAO, ST1,ST2, Z2}. 
Additional references regarding the $L_{p}$ Minkowski problem and Minkowski-type problems can be found, for example, in
\cite{BLYZ, WC, CW, GG, GL, GM, H1, HL1, HLYZ2, HLYZ1, HMS, JI, KL, LWA, LUT, LO, LYZ5, MIN, ST1, ST2, Z3, Zhu15}. Applications of the solutions to the $L_{p}$ Minkowski problem can be found in, e.g., \cite{AN1, AN2, KSC, Z, GH, LYZ4, CLYZ, HS1, HUS, Iva13, HS2, HSX, TUO}.\\

\noindent\textbf{$L_p$-Minkowski problem:} For $p\in\R$, what are the necessary and
sufficient conditions on a finite Borel measure $\mu$ on $S^{n-1}$ to ensure
that $\mu$ is the $L_p$ surface area measure of a convex body in
$\mathbb{R}^{n}$?\\

Besides discrete measures corresponding to polytopes, an important special case is when
\begin{equation}
\label{densityfunction}
d\mu=f\,d{\mathcal H}^{n-1}
\end{equation}
for some non-negative measurable function $f$ on $S^{n-1}$. 
If $p<1$ and (\ref{densityfunction}) holds, then the $L_p$-Minkowski problem amounts to solving the Monge-Amp\`ere type equation
\begin{equation}
\label{MongeAmper}
h^{1-p}\det(\nabla^2h+h I)=nf
\end{equation}
where $h$ is the unknown  non-negative function on $S^{n-1}$ to be found (the support function), $\nabla^2 h$ denotes the Hessian matrix of $h$ with respect to an orthonormal frame on $S^{n-1}$, and $I$ is the identity matrix.

If $n=2$, then we may assume that both $h$ and $f$ are non-negative periodic functions on $\R$ with period $2\pi$. In this case the corresponding differential equation is
\begin{equation}
\label{MongeAmperplanar}
h^{1-p}(h''+h)=2f.
\end{equation}
After earlier work by V. Umanskiy \cite{Uma03} and W. Chen \cite{WC}, equation (\ref{MongeAmperplanar}) 
in the $\pi$-periodic case that corresponds to planar origin symmetric convex bodies has been thoroughly investigated by M.Y. Jiang \cite{JI} if $p>-2$, and by M.N. Ivaki \cite{Iva13} if $p=-2$ (the "critical case"). 

Here we concentrate on the case $p\in(0,1)$. The case when $\mu$ has positive density function is handled by  Chou, Wang \cite{CW}:

\begin{theorem}[Chou, Wang]
\label{theodens}
If $p\in(-n,1)$,  $n\geq 2$ and $\mu$ is a Borel measure  on $S^{n-1}$  satisfying (\ref{densityfunction}) where $f$ is bounded and
$\inf_{u\in S^{n-1}} f(u)>0$, then $\mu$ 
is the $L_p$-surface area measure of a convex body  $K\in{\mathcal K}_{0}^n$.
\end{theorem}

We note that if $p\in(2-n,1)$,  then there exists $K\in{\mathcal K}_{0}^n$ with $o\in\partial K$ such that
$dS_{K,p}=f\,d{\mathcal H}^{n-1}$ for a positive continuous $f:\,S^{n-1}\to\R$ (see Example~\ref{oboundary}).

If $p\in(0,1)$, then the $L_p$-Minkowski problem for polytopes has been solved by Zhu \cite{Zhu15}. 

\begin{theorem}[Zhu]
\label{theopoly0<1}
For $p\in(0,1)$ and $n\geq 2$, a non-trivial discrete Borel measure $\mu$ on $S^{n-1}$ is the $L_p$-surface area measure of a polytope  $P\in{\mathcal K}_{0}^n$ 
with $o\in{\rm int}\,P$ if and only if $\mu$ is not concentrated on any closed hemisphere.
\end{theorem}
\noindent {\bf Remark } For the $\mu$ and $P$ as in Theorem~\ref{theopoly0<1}, if $G\subset O(n)$
 is a subgroup such that 
$\mu(\{Au\})=\mu(\{u\})$ for any $u\in S^{n-1}$ and $A\in G$, then 
 one may assume that $AP=P$ for any $A\in G$, as we explain in the Appendix.\\

For $p\in(0,1)$,  C. Haberl, E. Lutwak, D. Yang, G. Zhang \cite{HLYZ1} solved the $L_p$-Minkowski problem for even measures, or equivalently, for origin symmetric convex bodies.

\begin{theorem}[Haberl, Lutwak, Yang, Zhang]
\label{theoeven0<1}
For $p\in(0,1)$ and $n\geq 2$, a non-zero finite even Borel measure $\mu$ on $S^{n-1}$ is the $L_p$-surface area measure of an origin symmetric 
$K\in{\mathcal K}_{0}^n$ if and only if $\mu$ is not concentrated on any great subsphere.
\end{theorem}

The main goal of the paper is to solve the planar $L_p$ Minkowski problem in full generality if $p\in(0,1)$.
For a measure $\mu$ on $S^1$, we write ${\rm supp}\,\mu$ for its support. 

\begin{theorem}
\label{theoplanar0<1}
For $p\in(0,1)$ and a non-zero finite Borel measure $\mu$ on $S^1$, $\mu$ is the $L_p$-surface area measure of a convex body  $K\in{\mathcal K}_{0}^2$ if and only if ${\rm supp}\,\mu$ does not consist of a pair of antipodal vectors.
\end{theorem}
\noindent {\bf Remark } For the $\mu$ and $K$ as in Theorem~\ref{theoplanar0<1}, if $G\subset O(2)$ is a finite subgroup such that 
$\mu(A\omega)=\mu(\omega)$ for every Borel $\omega\subset S^1$ and $A\in G$, then one may assume that
$AK=K$ for any $A\in G$.\\

\begin{coro}
\label{coroplanar0<1}
For $p\in(0,1)$ and every non-negative $2\pi$-periodic function $f$  with
$0<\int_0^{2\pi}f<\infty$, the differential equation
(\ref{MongeAmperplanar}) has a non-negative $2\pi$-periodic weak solution.
\end{coro}
\noindent {\bf Remark } If the $f$ in (\ref{MongeAmperplanar}) is even, or is periodic with respect to $2\pi/k$ for an integer $k\geq 2$, then the solution $h$ can be also chosen even, or periodic with respect to $2\pi/k$, respectively.\\

Unfortunately, the method of the proof of Theorem~\ref{theoplanar0<1} does not extend to higher dimensions 
(see Example~\ref{highdimno}, and the remarks above). 

We note that for $p\in(2-n,1)$ in $\R^n$ (or $p\in(0,1)$ in $\R^2$), even if the function  $f$ on the right hand side of (\ref{MongeAmper}) or (\ref{MongeAmperplanar}) is positive and continuous, then possibly $o\in \partial K$ for the solution $K$. The following example is based on the example the end of Hug, Lutwak, Yang, Zhang \cite{HLYZ2}, and on examples in the preprint by Guan, Lin \cite{GL} and in Chou, Wang \cite{CW}.

\begin{example}
\label{oboundary}
If $p\in(2-n,1)$,  then there exists $K\in{\mathcal K}_{0}^n$ with $C^2$ boundary with $o\in\partial K$ such that
$dS_{K,p}=f\,d{\mathcal H}^{n-1}$ for a positive continuous $f:\,S^{n-1}\to\R$.
\end{example}
\proof
We fix $v\in S^{n-1}$, set
$B^{n-1}=v^\bot\cap B^n$ and for $x\in v^\bot$ and $t\in\R$, we write point $(x,t)=x+tv$.
For 
$$
q= \frac{2(n-1)}{n+p-2} >2,
$$
we consider the $C^2$ function $g(x)= \| x \|^q$ on $B^{n-1}$ where $\|\cdot\|$ stands for the Euclidean norm. We define the convex body $K$ 
in $\R^n$ with $C^2$ boundary in a way such that $o \in \partial K$ and the graph
$\{(x,g(x)):\,x\in B^{n-1}\}$ of $g$ above $B^{n-1}$
is a subset of $\partial K$. We may assume that $\partial K$ has positive Gau{\ss} curvature at each 
$z\in\partial K\backslash\{o\}$.

We observe that $K$ is strictly convex and $-v$ is the exterior unit normal at $o$, and hence $S_K(\{-v\})=0$.
If $z\in\partial K$, then we write $\nu(z)$ for the exterior unit normal at $z$,
and $\kappa(\nu(z))$ for the Gau{\ss} curvature at $z$, therefore even if $\kappa(-v)=0$, we have
$$
dS_K=\kappa^{-1}\,d\mathcal{H}^{n-1}.
$$
In turn, we deduce that
\begin{equation}
\label{SKpRadonNik}
dS_{K,p}=h_K^{1-p}\kappa^{-1}\,d\mathcal{H}^{n-1}.
\end{equation}

 Let $x\in B^{n-1}$ satisfy $0<\|x\|<1$, and let $z=(x,g(x))$, and hence $\kappa(\nu(z))>0$.
We have 
$$
\nabla g(x)  = q \|x \| ^{q-2} x\mbox{ \ and \ }
\nu(z)=a(x)^{-1}(\nabla g(x),-1)
$$
for
$$
a(x)= (1+\| \nabla g(x)\|^2)^{1/2}.	
$$		
In particular, writing $u=\nu(z)$, we have
$$
h_K(u)=\langle u,z\rangle=a(x)^{-1}\left(\langle \nabla g(x),x\rangle-g(x)\right)=
 a(x)^{-1}(q-1) \| x \|^q.
$$
In addition,
$$
\kappa(u)=a(x)^{-(n+1)} \det (\nabla ^2 g(x))= (q-1)q^{n-1} a(x)^{-(n+1)}  \| x \|^{(q-2)(n-1)},
$$ 
therefore  the Radon-Nikodym derivative in (\ref{SKpRadonNik}) is
$$ 
h_K(u)^{1-p}\kappa(u)^{-1} = 
(q-1)^{-p}q^{1-n}a(x)^{n+p} \|x \|^{q(1-p)-(q-2)(n-1)}=(q-1)^{-p}q^{1-n}a(x)^{n+p}. 
$$
Since $a(x)$ is continuous and positive function of  $x\in B^{n-1}$, we deduce that $S_{K,p}$ has a positive and continuous Radon-Nikodym derivative $f$ with respect to $\mathcal{H}^{n-1}$ on $S^{n-1}$. \ Q.E.D.\\

\noindent {\bf Note added in proof } In the meanwhile,
S. Chen, Q.-R. Li, G. Zhu \cite{CSL17} have essentially solved the $L_p$-Minkowski Problem for $0<p<1$ in all dimensions using a substantially different argument.

\section{Preliminary statements}

In this section, we prove some statements that are essential in proving Theorem~\ref{theoplanar0<1}.
For $v\in S^{n-1}$ and $t\in[0,1)$, let
$$
\Omega(v,t)=\{u \in S^{n-1}:\,\langle u,v\rangle>t\}.
$$
In particular, $\Omega(v,0)$ is the open hemi-sphere centered at $v$.

\begin{lemma}
\label{open-hemispheres}
If $\mu$ is a finite Borel measure on $S^{n-1}$ such that the measure of any open hemi-sphere is positive, then there exists 
$\delta\in(0,\frac12)$ such that for any $v\in S^{n-1}$, 
$$
\mu\left(\Omega(v,\delta)\right)>\delta.
$$
\end{lemma}
\noindent {\bf Remark } We may choose $\delta$ small enough such
that also $\mu(S^{n-1})<1/\delta$.\\
\proof
Suppose, to the contrary, that for any $k \in \mathbb{N},\,  k >1$, there exists $u_k \in S^{n-1}$ for which $\mu\left(\Omega\left( u_k,\frac{1}{k} \right) \right)\le \frac{1}{k}$. 
It follows from the compactness of $S^{n-1}$ that there is a convergent subsequence $\{ u_{k_j}\}$ of $\{u_k\}$ to some $u \in S^{n-1}$.  

Since $\mu\left(\Omega(u,0)\right) >0$, there exists $\tau= \cos \alpha$
for  $\alpha \in \left(  0, \frac{\pi}{2} \right) $ such that $\mu\left(\Omega(u,\tau)\right) >0$. There
exist large enough
 $k_j \in \mathbb{N}$ satisfying
$ \frac{1}{k_j} < \mu\left(\Omega(u,\tau)\right)$, $ \frac{1}{k_j} <\cos\frac{\pi+2\alpha}4$
and the angle $\theta$ of $u_{k_j}$ and $u$ is at most $\frac{\pi-2\alpha}4$. Since
$$
\cos(\alpha+\theta)\geq \cos\left(\alpha+\frac{\pi-2\alpha}4\right)=\cos\frac{\pi+2\alpha}4> \frac{1}{k_j},
$$
 the spherical triangle inequality yields $ \Omega(u,\tau)\subset\Omega \left( u_{k_j},\frac{1}{k_j} \right)$.
We deduce  that
$$
\mu\left(\Omega \left( u_{k_j},\frac{1}{k_j} \right)\right)\geq \mu\left(\Omega(u,\tau)\right)>\frac{1}{k_j} ,
$$
contradicting the definition of $u_k$, and proving Lemma~\ref{open-hemispheres}. Q.E.D.\\

Recall that the convex compact sets $K_m$ tend to the convex compact set $K$ in $\R^n$ (see R. Schneider \cite{SCH}) if 
$$
\lim_{m\to \infty}\max\{u\in S^{n-1}:\,\|h_{K_m}(u)-h_K(u)\|\}=0.
$$
We also note that the surface area measure can be extended to compact convex sets. Let $K$ be a compact convex set in $\R^n$. If ${\rm dim}\,K\leq n-2$, then $S_K$ is the constant zero measure. In addition,
if ${\rm dim}\,K=n-1$ and $v\in S^{n-1}$ is normal to ${\rm aff}\,K$,  then $S_K$ is concentrated on $\{\pm v\}$, and 
$S_K(\{v\})=S_K(\{-v\})=\mathcal{H}^{n-1}(K)$.

\begin{lemma}
\label{Lp-convergence}
If $\varphi:[0,\infty)\to[0,\infty)$ is continuous, and the sequence of convex compact convex sets $K_m$ with $o\in K_m$ tends to the convex compact set $K$ in $\R^n$, then the measures $\varphi\circ h_{K_m}\,dS_{K_m}$ tend weakly to $\varphi\circ h_{K}\,dS_{K}$.
\end{lemma}
\proof According to Theorem 4.2.1 in R. Schneider \cite{SCH}, $S_{K_m}$ tends weakly to $S_{K}$.
Since $o\in K_m$ for all $K_m$, we have $o\in K$. There exists $R>0$ such that $K_m\subset RB^n$ for every $m$, and hence $h_{K_m}(u)\leq R$ for $m$. Since $\varphi$ is uniformly continuous on $[0,R]$, for any continuous function 
$g:\,S^{n-1}\to \R$, the function $u\mapsto g(u)\varphi(h_{K_m}(u))$ tends uniformly
to $u\mapsto g(u)\varphi(h_{K}(u))$ on $S^{n-1}$. Therefore 
$(\varphi\circ h_{K_m})g\,dS_{K_m}$ tends  to $(\varphi\circ h_{K})g\,dS_{K}$. \ Q.E.D.\\

\begin{coro}
\label{Lp-convergence-cor}
If $p\leq 1$, and a sequence of compact convex sets $K_m$ with $o\in K_m$ tends to the compact convex set $K$ in $\R^n$,
then $S_{K_m,p}$ tends weakly to $S_{K,p}$.
\end{coro}

For $u_1,\ldots,u_k\in S^{n-1}$, we set 
$$
{\rm pos}\{u_1,\ldots,u_k\}=\{\lambda_1 u_1+\ldots+\lambda_k u_k:\,\lambda_1,\ldots,\lambda_k\geq 0\}.
$$

\begin{lemma}
\label{pos-hull}
If $x\in \R^n$ and $u_1,\ldots,u_k\in S^{n-1}$ satisfy that
$\langle u_i,x\rangle\geq 0$ for $i=1,\ldots,k$,  then for every $u\in S^{n-1}\cap {\rm pos}\{u_1,\ldots,u_k\}$, we have
$$
\langle u,x\rangle\geq \min\{\langle u_1,x\rangle,\ldots,\langle u_k,x\rangle\}.
$$ 
\end{lemma}
\proof We may assume that $\langle u_1,x\rangle\leq\langle u_i,x\rangle$  for $i=1,\ldots,k$.
The convexity of the unit ball yields that there exist  
$\lambda_1,\ldots,\lambda_k\geq 0$ with $\lambda_1+\ldots+\lambda_k\geq 1$ such that
$u=\lambda_1 u_1+\ldots+\lambda_k u_k$, and hence
$$
\langle u,x\rangle=\sum_{i=1}^k\lambda_i\langle u_i,x\rangle\geq 
\left(\sum_{i=1}^k\lambda_i\right)\langle u_1,x\rangle\geq \langle u_1,x\rangle.
\mbox{ \ \ \ Q.E.D.}
$$

For a planar convex body $K$ in $\R^2$, we say that $x_1,x_2\in\partial K$ are opposite points if there exists an exterior normal $u\in S^1$ at $x_1$ such that $-u$ is an exterior normal at $x_2\in \partial K$. If $x_1,x_2\in\partial K$ are not opposite, then we write $\sigma(K,x_1,x_2)$ for the arc of $\partial K$ connecting $x_1$ and $x_2$ not containing opposite points.  It is possible that $x_1=x_2$.
We observe that if $x\in \sigma(K,x_1,x_2)\backslash\{x_1,x_2\}$, then
\begin{equation}
\label{arcnormal}
\nu_K(x)\in {\rm pos}\{\nu_K(x_1),\nu_K(x_2)\}.
\end{equation}

\begin{claim}
\label{qzm}
For $p<1$, a planar convex body $K$ in $\R^2$ and non-opposite $x_1,x_2\in \partial K$,
if $\langle x_1,\nu_K(x_2)\rangle>0$ and $\langle x_2-x_1,u\rangle>0$ for $u\in S^1$, then
$$
\min\{h_K(\nu_K(x_1)),\langle x_1,\nu_K(x_2)\rangle\}^{1-p}\cdot  \langle x_2-x_1,u\rangle\leq 
\int_{S^1}h_K^{1-p}\,dS_K.
$$
\end{claim}
\proof
If $x\in \sigma(K,x_1,x_2)$ is a smooth point, then (\ref{arcnormal}) and 
 Lemma~\ref{pos-hull} yield
$$
\langle x,\nu_K(x)\rangle\geq  \langle x_1,\nu_K(x)\rangle\geq
\min\{\langle x_1,\nu_K(x_1)\rangle,\langle x_1, \nu_K(x_2)\rangle\}=
\min\{h_K(\nu_K(x_1)),\langle x_1,\nu_K(x_2)\rangle\}. 
$$
Therefore
\begin{eqnarray*}
\int_{S^1}h_K^{1-p}\,dS_K&=&\int_{\partial K}\langle x,\nu_{K}(x)\rangle^{1-p}\,d\mathcal{H}^1(x)>
\int_{\sigma(K,x_1,x_2)}\langle x,\nu_{K}(x)\rangle^{1-p}\,d\mathcal{H}^1(x) \\
&\geq&
\min\{h_K(\nu_K(x_1)),\langle x_1,\nu_K(x_2)\rangle\}^{1-p}\cdot 
\mathcal{H}^{1}(\sigma(K,x_1,x_2)),
\end{eqnarray*} 
and finally Claim~\ref{qzm} follows from $\mathcal{H}^{1}(\sigma(K,x_1,x_2))\geq \langle x_2-x_1,u\rangle$.
Q.E.D.

\section{Proof of Theorem~\ref{theoplanar0<1} if the measure of any open semicircle is positive} 
\label{secopensemipos}
Let $p\in(0,1)$, let $\mu$ be a finite Borel measure on $S^1$ such that the measure of any open semicircle is positive,
and let $\delta\in(0,\frac12)$ be the constant of  Lemma~\ref{open-hemispheres}  
for $\mu$ also satisfying $\mu(S^1)<1/\delta$.

We construct a sequence $\{\mu_m\}$ of discrete Borel measures on $S^1$ tending weakly to $\mu$ such that the $\mu_m$ measure of any open semicircle is positive for each $m$. It is the easiest to construct the sequence by identifying $\R^2$ with $\C$. For $m\geq 3$, we write $u_{jm}=e^{ij2\pi/m}$ for $i=\sqrt{-1}$ and $j=1,\ldots,m$, and we define
$\mu_m$ be the measure having the support $\{u_{1m},\ldots,u_{mm}\}$ with
$$
\mu_m(\{u_{jm}\})=\frac1{m^2}+
\mu\left(\{e^{it}:\,(j-1)2\pi<t\leq j2\pi\} \right)\mbox{ \ for $j=1,\ldots,m$}.
$$

According to Theorem~\ref{theopoly0<1} due to Zhu \cite{Zhu15}, there exists a polygon $P_m$ with 
$o\in{\rm int}\,P_m$ such that $d\mu_m=h_{P_m}^{1-p}\,dS_{P_m}$ for each $m$. It follows from 
 Lemma~\ref{Lp-convergence}
that we may assume that
\begin{equation}
\label{Pmtotal}
\int_{S^1}h_{P_m}^{1-p}\,dS_{P_m}<1/\delta.
\end{equation}

\begin{prop}
\label{Pmbounded}
$\{P_m\}$ is bounded.
\end{prop}
\proof
We assume that $d_m={\rm diam}\,P_m$ tends to infinity, and seek a contradiction. Choose $y_m,z_m\in P_m$ such that $\|z_m-y_m\|=d_m$ and $\|z_m\|\geq \|y_m\|$. Let $v_m=(z_m-y_m)/\|z_m-y_m\|$, and let $w_m\in S^1$ be orthogonal to $v_m$. Since $[y_m,z_m]$ is a diameter of $P_m$,  $v_m$ and $-v_m$ are exterior normals at $z_m$ and $y_m$, respectively. 
It follows that $\langle z_m, v_m\rangle\geq d_m/2$. By possibly taking subsequences, we may assume that $v_m$ tends to $\tilde{v}\in S^1$. It follows from Lemma~\ref{open-hemispheres} 
and Lemma~\ref{Lp-convergence} that if $m$ is large, then
\begin{equation}
\label{Pmvm}
\int_{\Omega(-v_m,\delta/2)}h_{P_m}^{1-p}\,dS_{P_m}>\delta/2.
\end{equation}
We prove Proposition~\ref{Pmbounded} based on the series of auxiliary statements Lemma~\ref{Lemma1} to Lemma~\ref{a*b*}.

 Let $a_m,b_m\in\partial P_m$ such that
$\langle a_m-b_m,w_m\rangle >0$ and $\langle a_m,v_m\rangle =\langle b_m,v_m\rangle =d_m/4$. 
We also deduce that $[a_m,b_m]\cap{\rm int}\,P_m\neq \emptyset$ for the segment $[a_m,b_m]$, and the definition of $a_m,b_m$,
$\langle z_m, v_m\rangle\geq d_m/2$ and 
$\langle y_m, -v_m\rangle\geq 0$ imply
$$
\mbox{$\langle z_m-a_m,v_m\rangle=\langle z_m-b_m,v_m\rangle\geq d_m/4$ and 
$\langle y_m-a_m,-v_m\rangle=\langle y_m-b_m,-v_m\rangle\geq d_m/4$}.
$$

\begin{lemma}
\label{Lemma1}
There exists $c_1>0$ depending on $\mu$ and $p$ such that if $m$ is large, then
\begin{equation}
\label{nupmqm}
h_{P_m}(\nu_{P_m}(a_m))\leq c_1d_m^{\frac{-1}{1-p}}
\mbox{ \ and \ } h_{P_m}(\nu_{P_m}(b_m))\leq c_1d_m^{\frac{-1}{1-p}}.
\end{equation}
\end{lemma}
\proof 
Since $\langle z_m-a_m,v_m\rangle\geq d_m/4$,
(\ref{Pmtotal}) and Claim~\ref{qzm} with $x_1=a_m$, $x_2=z_m$ and $v=v_m$ yield (\ref{nupmqm}) for $\nu_{P_m}(a_m)$, and the analogous argument works 
in the case of $\nu_{P_m}(b_m)$ with $a_m$ is replaced by $b_m$. Q.E.D.\\

Our intermediate goal, from Lemma~\ref{Lemma2} to Lemma~\ref{nuabmlow} is to
 show that $\nu_{P_m}(a_m)$ and $\nu_{P_m}(b_m)$ point essentially to the same direction as 
$w_m$ and $-w_m$, respectively, or in other words, 
$$
\lim_{m\to \infty}\langle \nu_{P_m}(a_m),v_m\rangle=\lim_{m\to \infty}\langle \nu_{P_m}(b_m),v_m\rangle=0.
$$
We frequently use the fact that
\begin{equation}
\label{x0x}
\langle \nu_{P_m}(x_0),x_0-x\rangle\geq 0
\end{equation}
for $x_0\in\partial P_m$ and $x\in P_m$.
In particular, $\langle \nu_{P_m}(a_m),w_m\rangle>0$ and $\langle \nu_{P_m}(b_m),-w_m\rangle>0$  as $\langle \nu_{P_m}(a_m),a_m-b_m\rangle> 0$
and $\langle \nu_{P_m}(b_m),b_m-a_m\rangle> 0$, respectively, by (\ref{x0x}) and $[a_m,b_m]\cap{\rm int}\,P_m\neq \emptyset$.

Below we frequently use the fact that if $p,q$ is an orthonormal basis for $\R^2$ and $x,y\in\R^2$, then
$$
\langle x,y\rangle=\langle x,p\rangle \langle y,p\rangle+\langle x,q\rangle \langle y,q\rangle.
$$

\begin{lemma}
\label{Lemma2}
For any $P_m$, we have
\begin{equation}
\label{nuamvmwm}
\frac{|\langle \nu_{P_m}(a_m),v_m\rangle|}{\langle \nu_{P_m}(a_m),w_m\rangle}\leq 
\frac{\langle a_m-z_m,w_m\rangle}{d_m/4}
\mbox{ \ and \ } \frac{|\langle \nu_{P_m}(b_m),v_m\rangle|}{\langle \nu_{P_m}(b_m),-w_m\rangle}\leq 
\frac{\langle b_m-z_m,-w_m\rangle}{d_m/4}.
\end{equation}
\end{lemma}
\proof 
It is enough to verify the statement about $\nu_{P_m}(a_m)$
where the definition of $a_m$  implies $\langle a_m-z_m,v_m\rangle\leq -d_m/4$. If $\langle \nu_{P_m}(a_m),v_m\rangle\geq 0$,
then 
\begin{eqnarray*}
0&\leq& \langle \nu_{P_m}(a_m),a_m-z_m\rangle=
\langle \nu_{P_m}(a_m),v_m\rangle\langle a_m-z_m,v_m\rangle
+\langle \nu_{P_m}(a_m),w_m\rangle\langle a_m-z_m,w_m\rangle\\
&\leq&
-\langle \nu_{P_m}(a_m),v_m\rangle(d_m/4)
+\langle \nu_{P_m}(a_m),w_m\rangle\langle a_m-z_m,w_m\rangle
\end{eqnarray*}
yields (\ref{nuamvmwm}). If $\langle \nu_{P_m}(a_m),-v_m\rangle\geq 0$,
then using $\langle a_m-y_m,-v_m\rangle\leq -d_m/4$ and
$\langle a_m-y_m,w_m\rangle=\langle a_m-z_m,w_m\rangle$, we deduce
\begin{eqnarray*}
0&\leq& \langle \nu_{P_m}(a_m),a_m-y_m\rangle=
\langle \nu_{P_m}(a_m),-v_m\rangle\langle a_m-y_m,-v_m\rangle
+\langle \nu_{P_m}(a_m),w_m\rangle\langle a_m-y_m,w_m\rangle\\
&\leq&
-\langle \nu_{P_m}(a_m),-v_m\rangle(d_m/4)
+\langle \nu_{P_m}(a_m),w_m\rangle\langle a_m-y_m,w_m\rangle,
\end{eqnarray*}
and in turn we have (\ref{nuamvmwm}). \ Q.E.D.\\

\begin{coro}
For any $P_m$, we have
\begin{equation}
\label{nuamwmest}
\langle \nu_{P_m}(a_m),w_m\rangle\geq \frac15
\mbox{ \ and \ } \langle \nu_{P_m}(b_m),-w_m\rangle\geq \frac15.
\end{equation}
\end{coro}
\proof It is enough to verify the statement about $\nu_{P_m}(a_m)$. Let $\gamma_m=\angle(\nu_{P_m}(a_m),w_m)$.
Since $\langle a_m-z_m,w_m\rangle\leq d_m$ and $\langle b_m-z_m,-w_m\rangle\leq d_m$ follow from
$\|a_m-b_m\|\leq d_m$, we conclude from (\ref{nuamvmwm}) that $\tan\gamma_m\leq 4$. We deduce that
$$
\langle \nu_{P_m}(a_m),w_m\rangle=\cos\gamma_m=(1+\tan^2\gamma_m)^{-1/2}\geq \frac1{\sqrt{17}}>\frac15.
\mbox{ \ Q.E.D.}\\
$$

Possibly interchanging $w_m$ with $-w_m$, and the role of $a_m$ and $b_m$, we may assume that
$\langle y_m,w_m\rangle=\langle z_m,w_m\rangle\geq 0$. We have $z_m=t_mv_m+r_mw_m$ and 
$y_m=s_m(-v_m)+r_mw_m$ for $t_m\geq s_m\geq 0$ and $r_m\geq 0$ where
$t_m+s_m=d_m$. In particular, $t_m\geq d_m/2$.

\begin{lemma}
\label{Lemma3}
There exist $c_2,c_3,c_4>0$ depending on $\mu$ and $p$ such that if $m$ is large, then
\begin{eqnarray}
\label{rm}
r_m &\leq &c_2 d_m^{\frac{-1}{1-p}}\\
\label{nuamupp}
\langle v_m,\nu_{P_m}(a_m)\rangle &\leq &c_3 d_m^{\frac{p-2}{1-p}}\\
\label{nubmupp}
\langle v_m,\nu_{P_m}(b_m)\rangle &\leq &c_4 d_m^{\frac{p-2}{1-p}}.
\end{eqnarray}
\end{lemma}
\proof
If $\langle v_m,\nu_{P_m}(a_m)\rangle\geq 0$, then (\ref{nupmqm})  implies
$$
r_m\langle w_m,\nu_{P_m}(a_m)\rangle+t_m\langle v_m,\nu_{P_m}(a_m)\rangle
= \langle z_m,\nu_{P_m}(a_m)\rangle \leq \langle a_m,\nu_{P_m}(a_m)\rangle
\leq c_1 d_m^{\frac{-1}{1-p}},
$$
which in turn yields (\ref{rm}) by (\ref{nuamwmest}) in this case,
and in addition, yields (\ref{nuamupp}) by $t_m\geq d_m/2$. Similarly, if $\langle -v_m,\nu(a_m)\rangle\geq 0$, then we have
$$
r_m\langle w_m,\nu_{P_m}(a_m)\rangle+s_m\langle -v_m,\nu_{P_m}(a_m)\rangle
= \langle y_m,\nu_{P_m}(a_m)\rangle \leq \langle a_m,\nu_{P_m}(a_m)\rangle
\leq c_1 d_m^{\frac{-1}{1-p}},
$$
and we conclude (\ref{rm}) using again (\ref{nuamwmest}).  Finally, if 
$\langle v_m,\nu_{P_m}(b_m)\rangle\geq 0$, then combining $0<\langle -w_m,\nu_{P_m}(b_m)\rangle\leq 1$, (\ref{rm}) and
$$
-r_m\langle -w_m,\nu_{P_m}(b_m)\rangle+t_m\langle v_m,\nu_{P_m}(b_m)\rangle
= \langle z_m,\nu_{P_m}(b_m)\rangle \leq \langle b_m,\nu_{P_m}(b_m)\rangle
\leq c_1 d_m^{\frac{-1}{1-p}}
$$
implies
$$
t_m\langle v_m,\nu_{P_m}(b_m)\rangle \leq (c_1+c_2) d_m^{\frac{-1}{1-p}},
$$
and in turn we conclude (\ref{nubmupp}) by $t_m\geq d_m/2$. Q.E.D.\\

\begin{lemma}
\label{nuabmlow}
There exist $c_5,c_6>0$ depending on $\mu$ and $p$ such that
 if $m$ is large, then
\begin{eqnarray}
\label{nuamlow}
\langle v_m,\nu_{P_m}(a_m)\rangle &\geq &-c_5 d_m^{\frac{p-1}{3-3p+p^2}-1}\\
\label{nubmlow}
\langle v_m,\nu_{P_m}(b_m)\rangle &\geq &-c_6 d_m^{\frac{p-1}{3-3p+p^2}-1}.
\end{eqnarray}
\end{lemma}
\proof According to (\ref{nuamvmwm}), it is sufficient to prove that
there exist $c_7,c_8>0$ depending on $\mu$ and $p$ such that
\begin{eqnarray}
\label{nuamheightupp}
\alpha_m=\langle a_m-z_m,w_m\rangle &\leq &c_7 d_m^{\frac{p-1}{3-3p+p^2}}
\mbox{ \ provided $\langle v_m,\nu_{P_m}(a_m)\rangle<0$,}\\
\label{nubmheightupp}
\beta_m=\langle b_m-z_m,-w_m\rangle &\leq &c_8 d_m^{\frac{p-1}{3-3p+p^2}}
\mbox{ \ provided $\langle v_m,\nu_{P_m}(b_m)\rangle<0$.}
\end{eqnarray}
For (\ref{nuamheightupp}),  $\|a_m-z_m\|\leq d_m$ and $|\langle a_m-z_m,v_m\rangle|\geq d_m/4$ yield
$\alpha_m\leq\frac{\sqrt{15}}4\,d_m$, and hence 
$$
\eta_m=\left( \frac{\alpha_m}{d_m}\right)^{\frac{1-p}{2-p}}\leq
\left( \frac{\sqrt{15}}{4}\right)^{\frac{1-p}{2-p}}<1.
$$
The constant $\eta_m$ is chosen in a way such that the calculations in Case 1 and in Case 2 lead to the same estimate up to a constant factor.

We consider the vector $e_m\in S^1$ such that $\langle e_m, v_m\rangle=\eta_m$ and 
$\langle e_m, w_m\rangle>0$, and hence
there exists $c_9>0$ depending on  $p$ such that
$$
 \langle  e_m, w_m\rangle\geq c_9.
$$
There exists $a'_m\in\sigma(P_m,a_m,z_m)$ such that $w_m$ is an exterior unit normal, and there exists
 $\tilde{a}_m\in\sigma(P_m,a'_m,z_m)$ such that 
$e_m$ is an exterior unit normal at $\tilde{a}_m$. In particular, we may assume that $\nu_{P_m}(a'_m)=w_m$ and $\nu_{P_m}(\tilde{a}_m)=e_m$, and we have
$$
\langle a'_m,w_m\rangle \geq \langle a'_m-z_m,w_m\rangle=h_{P_m}(w_m)-\langle z_m,w_m\rangle\geq 
\langle a_m-z_m,w_m\rangle=\alpha_m.
$$
We distinguish two cases. \\

\noindent{\bf Case 1} $\langle \tilde{a}_m-z_m,w_m\rangle< \alpha_m/2$

We want to apply
Claim~\ref{qzm} with $x_1=a'_m$ $x_2=\tilde{a}_m$ and $u=v_m$. 
Since both of $\langle a'_m,w_m\rangle$ and $\langle e_m,w_m\rangle$
are positive, and  $\langle a'_m,v_m\rangle\geq d_m/4$, 
$\langle e_m, v_m\rangle=\eta_m$ 
and $d_m\geq \alpha_m$, we deduce that
$$
\langle a'_m,e_m\rangle=
\langle a'_m,v_m\rangle\langle e_m, v_m\rangle+\langle a'_m,w_m\rangle\langle e_m, w_m\rangle
\geq 
(d_m/4)\eta_m=\mbox{$\frac14$}\,\alpha_m^{\frac{1-p}{2-p}}d_m^{\frac{1}{2-p}} \geq \frac{\alpha_m}4.
$$
In addition, $h_{P_m}(w_m)\geq \alpha_m$, thus 
$\min\{h_{P_m}(w_m),\langle a'_m,e_m\rangle\}\geq \frac{\alpha_m}4$.
 Since $\langle \tilde{a}_m-a'_m,w_m\rangle<-\alpha_m/2$ by the condition in Case 1, we have
$$
0\leq \langle \tilde{a}_m-a'_m,e_m\rangle=
\langle \tilde{a}_m-a'_m,v_m\rangle\langle e_m,v_m\rangle+\langle \tilde{a}_m-a'_m,w_m\rangle\langle e_m,w_m\rangle\leq
\langle \tilde{a}_m-a'_m,v_m\rangle\eta_m-\frac{c_9\alpha_m}2,
$$
and hence
$$
\langle \tilde{a}_m-a'_m,v_m\rangle\geq \frac{c_9\alpha_m}{2\eta_m}=\frac{c_9}{2}\,\alpha_m^{\frac{1}{2-p}}d_m^{\frac{1-p}{2-p}}.
$$
Therefore  (\ref{Pmtotal}) and Claim~\ref{qzm} with $x_1=a'_m$, $x_2=\tilde{a}_m$ and $u=v_m$ imply
$$
\left( \frac{\alpha_m}4\right)^{1-p}\cdot  \frac{c_9}{2}\,\alpha_m^{\frac{1}{2-p}}d_m^{\frac{1-p}{2-p}} <\frac1{\delta},
$$
and in turn we conclude (\ref{nuamheightupp}). \\

\noindent{\bf Case 2} $\langle \tilde{a}_m-z_m,w_m\rangle\geq \alpha_m/2$

Now $\langle z_m,e_m\rangle\geq (d_m/4)\eta_m$ by $\langle z_m,w_m\rangle\geq 0$, thus 
$h_{P_m}(v_m)\geq d_m/2$ yields 
$$
\min\{h_{P_m}(v_m),\langle z_m,e_m\rangle\}\geq (d_m/4)\eta_m=
\frac14\,\alpha_m^{\frac{1-p}{2-p}}d_m^{\frac{1}{2-p}}.
$$ 
Therefore  (\ref{Pmtotal}) and Claim~\ref{qzm} with $x_1=z_m$, $x_2=\tilde{a}_m$ and $u=w_m$ yield
$$
\left(\frac14\,\alpha_m^{\frac{1-p}{2-p}}d_m^{\frac{1}{2-p}}\right)^{1-p}\cdot\frac{\alpha_m}2 <\frac1{\delta},
$$
and  we finally conclude (\ref{nuamheightupp}).

Next we turn to (\ref{nubmheightupp}) where 
the argument is similar to the argument for  (\ref{nuamheightupp}).  
The difference between the proofs of (\ref{nubmheightupp}) and (\ref{nuamheightupp}) is that now
 $\langle z_m,-w_m\rangle <0$. However, 
$\langle z_m,-w_m\rangle=-r_m>-c_2 d_m^{\frac{-1}{1-p}}$ according to (\ref{rm}).  
If 
$$
\beta_m< d_m^{\frac{p-1}{3-3p+p^2}},
$$
then (\ref{nubmheightupp}) readily holds. Therefore, we assume that
$$
\beta_m\geq d_m^{\frac{p-1}{3-3p+p^2}}.
$$
Since $\frac{-1}{1-p}<\frac{p-1}{3-3p+p^2}$, we may assume that $m$ is large enough to ensure that
\begin{equation}
\label{betamlarge0}
\beta_m\geq d_m^{\frac{p-1}{3-3p+p^2}}>4c_2 d_m^{\frac{-1}{1-p}}\geq 4r_m.
\end{equation}
 In particular, if $m$ is large, then
\begin{equation}
\label{betamlarge}
 \langle b_m,-w_m\rangle\geq \frac{3\beta_m}4.
\end{equation}

Since $\|b_m-z_m\|\leq d_m$ and $|\langle b_m-z_m,v_m\rangle|\geq d_m/4$ yield
$\beta_m\leq\frac{\sqrt{15}}4\,d_m$, we have
$$
\theta_m=\left( \frac{\beta_m}{d_m}\right)^{\frac{1-p}{2-p}}\leq
\left( \frac{\sqrt{15}}{4}\right)^{\frac{1-p}{2-p}}<1.
$$
We consider the vector $f_m\in S^1$ such that $\langle f_m, v_m\rangle=\theta_m$ and 
$\langle f_m, -w_m\rangle>0$, and hence for the
 $c_9>0$ above depending on  $p$, we have
$$
 \langle  f_m, -w_m\rangle\geq c_9.
$$
There exists $b'_m\in\sigma(P_m,b_m,z_m)$ such that $-w_m$ is an exterior unit normal, and there exists
 $\tilde{b}_m\in\sigma(P_m,b'_m,z_m)$ such that 
$f_m$ is an exterior unit normal at $\tilde{b}_m$. In particular, we may assume that $\nu_{P_m}(b'_m)=-w_m$ and $\nu_{P_m}(\tilde{b}_m)=f_m$, and we have
$$
 \langle b'_m-z_m,-w_m\rangle\geq \beta_m.
$$

Again, we distinguish two cases. \\

\noindent{\bf Case 1'} $\langle \tilde{b}_m-z_m,-w_m\rangle< \beta_m/2$

In this case, we are going to apply Claim~\ref{qzm} with $x_1=b'_m$, $x_2=\tilde{b}_m$ and $u=v_m$.
Since both of $\langle b'_m,-w_m\rangle$ and $\langle f_m,-w_m\rangle$
are positive, and  $\langle b'_m,v_m\rangle\geq d_m/4$, 
$\langle f_m, v_m\rangle=\theta_m$ 
and $d_m\geq \beta_m$, we deduce that
$$
\langle b'_m,f_m\rangle=
\langle b'_m,v_m\rangle\langle f_m, v_m\rangle+\langle b'_m,-w_m\rangle\langle f_m, -w_m\rangle
\geq 
(d_m/4)\theta_m=\mbox{$\frac14$}\,\beta_m^{\frac{1-p}{2-p}}d_m^{\frac{1}{2-p}} \geq \frac{\beta_m}4.
$$
In addition, $h_{P_m}(-w_m)\geq 3\beta_m/4$ by (\ref{betamlarge}), thus 
$\min\{h_{P_m}(-w_m),\langle b'_m,f_m\rangle\}\geq \frac{\beta_m}4$.
 Since 
$\langle \tilde{b}_m-b'_m,-w_m\rangle<-\beta_m/2$ by the condition in Case 1', we have
$$
0\leq \langle \tilde{b}_m-b'_m,f_m\rangle=
\langle \tilde{b}_m-b'_m,v_m\rangle\langle f_m,v_m\rangle+\langle \tilde{b}_m-b'_m,-w_m\rangle
\langle f_m,-w_m\rangle\leq
\langle \tilde{b}_m-b'_m,v_m\rangle\theta_m-\frac{c_9\beta_m}2,
$$
and hence
$$
\langle \tilde{b}_m-b'_m,v_m\rangle\geq \frac{c_9\beta_m}{2\theta_m}=
\frac{c_9}{2}\,\beta_m^{\frac{1}{2-p}}d_m^{\frac{1-p}{2-p}}.
$$
Therefore  (\ref{Pmtotal}) and Claim~\ref{qzm} with $x_1=b'_m$, $x_2=\tilde{b}_m$ and $u=v_m$ imply
$$
\left( \frac{\beta_m}4\right)^{1-p}\cdot  \frac{c_9}{2}\,\beta_m^{\frac{1}{2-p}}d_m^{\frac{1-p}{2-p}} <\frac1{\delta},
$$
and in turn we conclude (\ref{nubmheightupp}). \\

\noindent{\bf Case 2'} $\langle \tilde{b}_m-z_m,-w_m\rangle\geq \beta_m/2$

In this case, (\ref{rm}) implies
$$
\langle z_m,f_m\rangle=\langle z_m,v_m\rangle\langle f_m,v_m\rangle+
\langle z_m,-w_m\rangle\langle f_m,-w_m\rangle\geq (d_m/4)\theta_m
-c_2 d_m^{\frac{-1}{1-p}}.
$$  
Here, if $m$ is large, then using the definition of $\theta_m$ and (\ref{betamlarge0}), we have
$$
d_m\theta_m=d_m\left( \frac{\beta_m}{d_m}\right)^{\frac{1-p}{2-p}}=\beta_m^{\frac{1-p}{2-p}}d_m^{\frac{1}{2-p}}
\geq \left(4c_2 d_m^{\frac{-1}{1-p}}\right)^{\frac{1-p}{2-p}}d_m^{\frac{1}{2-p}}=(4c_2)^{\frac{1-p}{2-p}}
>8c_2 d_m^{\frac{-1}{1-p}},
$$
thus  $\langle z_m,f_m\rangle\geq (d_m/8)\eta_m$. It follows from $h_{P_m}(v_m)\geq d_m/2$ that
$$
\min\{h_{P_m}(v_m),\langle z_m,f_m\rangle\}\geq (d_m/8)\theta_m=
\frac18\,\beta_m^{\frac{1-p}{2-p}}d_m^{\frac{1}{2-p}}.
$$ 
Therefore  (\ref{Pmtotal}) and Claim~\ref{qzm} with $x_1=z_m$, $x_2=\tilde{b}_m$ and $u=-w_m$ yield
$$
\left(\frac18\,\beta_m^{\frac{1-p}{2-p}}d_m^{\frac{1}{2-p}}\right)^{1-p}\cdot\frac{\beta_m}2 <\frac1{\delta},
$$
and  we finally conclude (\ref{nubmheightupp}), and in turn Lemma~\ref{nuabmlow}. Q.E.D.\\

In order to finish the proof of Proposition~\ref{Pmbounded}, let $a^*_m\in\partial P_m$ maximize 
$\langle a^*_m,w_m\rangle$ under the condition that the $\gamma_m\in S^1$ with 
$\langle \gamma_m,-v_m\rangle =\delta/2$ and $\langle \gamma_m,w_m\rangle>0$ is an exterior unit normal at $a^*_m$,
and let $b^*_m\in\partial P_m$ maximize $\langle b^*_m,-w_m\rangle$ under the condition that the 
$\xi_m\in S^1$ with 
$\langle \xi_m,-v_m\rangle =\delta/2$ and $\langle \xi_m,-w_m\rangle>0$ is an exterior unit normal at $b^*_m$.

\begin{lemma} 
\label{a*b*}
There exist $c_{10},c_{11}>0$ depending on $\mu$ and $p$ such that if $m$ is large, then
\begin{eqnarray}
\label{a*}
\langle a^*_m-y_m,w_m\rangle &\leq& c_{10}d_m^{\frac{-1}{1-p}} \\
\label{b*}
\langle b^*_m-y_m,-w_m\rangle &\leq & c_{11}d_m^{\frac{-1}{1-p}}. 
\end{eqnarray}
\end{lemma}
\proof For (\ref{a*}), $\langle y_m,w_m\rangle=r_m\geq 0$ yields
\begin{equation}
\label{a*yma*}
\langle a^*_m-y_m,w_m\rangle\leq \langle a^*_m,w_m\rangle.
\end{equation}
Since $\langle a^*_m-y_m,v_m\rangle\geq 0$ and $\langle a^*_m-y_m,w_m\rangle\geq 0$, we have
$$
0\leq \langle a^*_m-y_m,\gamma_m\rangle=\langle a^*_m-y_m,v_m\rangle\langle \gamma_m,v_m\rangle
+\langle a^*_m-y_m,w_m\rangle\langle \gamma_m,w_m\rangle\leq 
\frac{-\delta}2\langle a^*_m-y_m,v_m\rangle+\langle a^*_m-y_m,w_m\rangle.
$$
In turn (\ref{a*yma*}) implies
\begin{equation}
\label{a*ym}
\langle a^*_m-y_m,v_m\rangle\leq \frac2{\delta}\langle a^*_m-y_m,w_m\rangle\leq \frac2{\delta}\langle a^*_m,w_m\rangle.
\end{equation}

It follows from (\ref{nupmqm}) and (\ref{nuamwmest}) that
\begin{eqnarray}
\nonumber
c_1d_m^{\frac{-1}{1-p}}&\geq& h_{P_m}(\nu_{P_m}(a_m))\geq \langle a^*_m,\nu_{P_m}(a_m)\rangle=
\langle a^*_m,v_m\rangle\langle \nu_{P_m}(a_m),v_m\rangle+
\langle a^*_m,w_m\rangle\langle \nu_{P_m}(a_m),w_m\rangle\\
\label{a*support}
&\geq &
\langle a^*_m,v_m\rangle\langle \nu_{P_m}(a_m),v_m\rangle+
\langle a^*_m,w_m\rangle/5.
\end{eqnarray}
The rest of the argument is divided into three cases according to the signs of 
$\langle a^*_m,v_m\rangle$ and $\langle \nu_{P_m}(a_m),v_m\rangle$.\\

\noindent{\bf Case 1 } $\langle a^*_m,v_m\rangle\langle \nu_{P_m}(a_m),v_m\rangle\geq 0$ 

In this case (\ref{a*support}) 
and  (\ref{a*yma*}) yield  (\ref{a*}) directly.\\

\noindent{\bf Case 2 } $\langle a^*_m,v_m\rangle>0$ and $\langle \nu_{P_m}(a_m),v_m\rangle<0$

In this case $\langle y_m,v_m\rangle\leq 0$
and (\ref{a*ym}) imply $\langle a^*_m,v_m\rangle\leq \frac2{\delta}\langle a^*_m,w_m\rangle$.
Since $|\langle \nu_{P_m}(a_m),v_m\rangle|<\frac{\delta}{20}$ for large $m$ according to
(\ref{nuamlow}), we conclude from (\ref{a*support}) that
$$
c_1d_m^{\frac{-1}{1-p}}\geq -\frac2{\delta}\langle a^*_m,w_m\rangle\cdot \frac{\delta}{20}
+\frac{\langle a^*_m,w_m\rangle}{5}=\frac{\langle a^*_m,w_m\rangle}{10},
$$
proving (\ref{a*}) by (\ref{a*yma*}). \\

\noindent{\bf Case 3 } $\langle a^*_m,v_m\rangle<0$ and 
$\langle \nu_{P_m}(a_m),v_m\rangle>0$

In this case, we have $\langle a^*_m,v_m\rangle\geq -d_m$ on the one hand, and
(\ref{nuamupp}) implies $\langle \nu_{P_m}(a_m),v_m\rangle<c_3 d_m^{\frac{p-2}{1-p}}$
on the other hand, therefore (\ref{a*support}) yields
$$
c_1d_m^{\frac{-1}{1-p}}\geq -d_m c_3 d_m^{\frac{p-2}{1-p}}+\frac{\langle a^*_m,w_m\rangle}{5}
= -c_3d_m^{\frac{-1}{1-p}}+\frac{\langle a^*_m,w_m\rangle}{5},
$$
completing the proof of (\ref{a*}) by (\ref{a*yma*}).\\

For (\ref{b*}), we may assume that
$$
\langle b^*_m-y_m,-w_m\rangle \geq 2 c_2d_m^{\frac{-1}{1-p}},
$$
otherwise (\ref{b*}) readily holds with $c_{11}=2c_2$.
Since $\langle b^*_m-y_m,-w_m\rangle\leq \langle b^*_m,-w_m\rangle+c_2d_m^{\frac{-1}{1-p}}$ by
(\ref{rm}), we have
\begin{equation}
\label{b*ym}
\langle b^*_m-y_m,-w_m\rangle\leq 2\langle b^*_m,-w_m\rangle.
\end{equation}
Therefore using (\ref{b*ym}) in place of $\langle a^*_m-y_m,w_m\rangle\leq \langle a^*_m,w_m\rangle$,
(\ref{b*}) can be proved similarly to (\ref{a*}), completing the proof of 
Lemma~\ref{a*b*}. Q.E.D.\\

Finally, to prove Proposition~\ref{Pmbounded}, we observe that combining Lemma~\ref{a*b*} with the definition of $a^*_m$ and $b^*_m$ yields that if $m$ is large, then
$$
\mathcal{H}^1(\sigma(P_m,a^*_m,b^*_m))\leq \frac2{\delta}\langle a^*_m-b^*_m,w_m\rangle\leq
\frac{2(c_{10}+c_{11})}{\delta}\cdot d_m^{\frac{-1}{1-p}}.
$$
It follows from applying first (\ref{Pmvm}), then $\langle x,\nu_{P_m}(x)\rangle\leq d_m$
for $x\in\partial K$  that if  $m$ is large, then
\begin{eqnarray*}
\frac{\delta}2&<&\int_{\Omega(-v_m,\delta/2)}h_{P_m}^{1-p}\,dS_{P_m}=
\int_{\sigma(P_m,a^*_m,b^*_m)}\langle x,\nu_{P_m}(x)\rangle^{1-p}\,d\mathcal{H}^1(x)\\
&\leq &d_m^{1-p}\cdot\mathcal{H}^1(\sigma(P_m,a^*_m,b^*_m))
\leq d_m^{1-p}\cdot \frac{2(c_{10}+c_{11})}{\delta}\cdot d_m^{\frac{-1}{1-p}}= \frac{2(c_{10}+c_{11})}{\delta}\cdot d_m^{\frac{p(p-2)}{1-p}},
\end{eqnarray*}
which is absurd as $\frac{p(p-2)}{1-p}<0$ and $d_m$ tends to infinity. This contradiction  verifies Proposition~\ref{Pmbounded}. \ Q.E.D.\\

\noindent{\bf Proof of Theorem~\ref{theoplanar0<1} if the measure of any open semicircle is positive } 
Since $\{P_m\}$ is bounded  and each $P_m$ contains the origin according to
Proposition~\ref{Pmbounded}, the Blaschke selection theorem provides a  subsequence $\{P_{m'}\}$ tending to a compact convex set $K$ with $o\in K$. It follows from Corollary~\ref{Lp-convergence-cor} that $S_{P_{m'},p}$ tends weakly to $S_{K,p}$. However, $\mu_{m'}=S_{P_{m'},p}$ tends weakly to $\mu$ by construction. Therefore 
$\mu=S_{K,p}$. 
Since any open semi-circle of $S^1$ has positive $\mu$ measure, we conclude that ${\rm int}\,K\neq \emptyset$. 

Finally, to prove the Remark after Theorem~\ref{theoplanar0<1},  let $G\subset O(2)$ be a finite subgroup such that 
$\mu(A\omega)=\mu(\omega)$ for any Borel $\omega\subset S^1$ and $A\in G$. The idea is that for large $m$, we subdivide $S^1$ into arcs of length less than $2\pi/m$ in a way such that the subdivision is symmetric with respect to $G$ and each endpoint has $\mu$ measure $0$.

We fix a regular $l$-gon $Q$, 
$l\geq 3$ whose vertices lie on $S^1$ such that $G$ is a subgroup of the symmetry group of $Q$. 
In addition, we consider the set $\Sigma$ of atoms of $\mu$; namely, 
the set of all $u\in S^1$ such that $\mu(\{u\})>0$. In particular, $\Sigma$ is countable. 

For $m\geq 2$, let $Q_m$ be a regular polygon with $lm$ vertices such that all vertices of $Q$ are vertices of $Q_m$, and let $G_m$ be the symmetry group of $Q_m$. We observe that $G_m$ contains rotations by angle $\frac{2\pi}{lm}$. We write $\Sigma_m$ for the set obtained from repeated applications of the elements of $G_m$  to the elements of $\Sigma$, and hence $\Sigma_m$ is countable, as well. For a fixed $x_0\in S^1\backslash \Sigma_m$, we consider the orbit $G_mx_0=\{Ax_0:\,A\in G_m\}$, and let $\mathcal{I}_m$
be the set of open arcs of $S^1$ that are the components of $S^1\backslash G_mx_0$. We observe that
$G_mx_0$ is disjoint from $\Sigma_m$, and hence $\mu(\sigma)=\mu({\rm cl}\,\sigma)$ for
$\sigma\in \mathcal{I}_m$. 

Now we define $\mu_m$. It is concentrated on the set of midpoints of all $\sigma\in \mathcal{I}_m$, and the $\mu_m$ measure of the midpoints of a $\sigma\in \mathcal{I}_m$ is $\mu(\sigma)$. In particular, 
$\mu_m$ is invariant under $G_m$, and hence $\mu_m$ is invariant under $G$. Since the length of each arc in 
$\mathcal{I}_m$ is at most $\frac{2\pi}{lm}$, we deduce that $\mu_u$ tends weakly to $\mu$.

According to the Remark after Theorem~\ref{theopoly0<1} due to Zhu \cite{Zhu15}, we may assume that each $P_m$ is invariant under $G$. Now the argument above shows that some subsequence of $\{P_m\}$ tends to a convex body $K$ satisfying $S_{K,p}=\mu$, and readily $K$  is invariant under $G$. \ Q.E.D.\\

Unfortunately, the proof of Theorem~\ref{theoplanar0<1} we present does not extend to higher dimensions. What we actually prove in this section (see Proposition~\ref{Pmbounded}) is the following statement: If $0<p<1$, $\mu$ is a bounded Borel measures on $S^1$ such that the $\mu$ measure of any open semi-circle is positive, and $P_m\in\mathcal{K}^2_o$ is a sequence of convex bodies such that $S_{P_m,p}$ tends weakly to $\mu$,
then the sequence $\{P_m\}$ is bounded. The following Example~\ref{highdimno} shows that this statement already fails in $n=3$ dimension. For $x_1,\ldots,x_k\in \R^3$, we write $[x_1,\ldots,x_k]$ for their convex hull.

\begin{example}
\label{highdimno}
For $p\in(0,1)$, there exist a  measure $\mu$ on $S^2$ such that any open hemisphere has positive measure and an unbounded sequence of polytopes $\{P_m\}$ in $\R^3$ such that
$o\in{\rm int} P_m$ and $S_{P_m,p}$ tends weakly to $\mu$.
\end{example}
\proof We define 
$$
u_0=(1,0,0),\;u_1=\left(\frac{-1}{\sqrt{2}},\frac{1}{\sqrt{2}},0\right),\;
u_2=\left(\frac{-1}{\sqrt{2}},\frac{-1}{\sqrt{2}},0\right),\;
u^+=(0,0,1),\;u^-=(0,0,-1),
$$
and the discrete measure $\mu$ with
 ${\rm supp}\,\mu=\{u_0,u_1,u_2,u^+,u^-\}$ and
$$
\mu(\{u_0\})=8,\;\mu(\{u_1\})=\mu(\{u_2\})=2^{\frac{p}2},\;\mu(\{u^+\})=\mu(\{u^-\})=3,
$$
and hence any open hemisphere has positive measure.

For $m\geq 2$ and $a=a_m=m^{-(2-p)}$, let
$$
v_{1,m}=(0,m,0),\;v^+_{1,m}=(m,2m,a),\;v^-_{1,m}=(m,2m,-a),
$$
$$
v_{2,m}=(0,-m,0),\;v^+_{2,m}=(m,-2m,a),\;v^-_{2,m}=(m,-2m,-a),
$$
and let $\widetilde{P}_m$ be their convex hull.
The exterior unit normals of the faces $F_{0,m}=[v^+_{i,m},v^-_{i,m}]_{i=1,2}$,
$F_{1,m}=[v_{1,m},v^+_{1,m},v^-_{1,m}]$ and
$F_{2,m}=[v_{2,m},v^+_{2,m},v^-_{2,m}]$  
 are $u_0,u_1,u_2$,
 respectively, which vectors are independent of $m$. In addition, $\widetilde{P}_m$
has two more facets, $F_{m}^+=[v_{1,m},v_{2,m},v^+_{1,m},v^+_{2,m}]$ and
$F_{m}^-=[v_{1,m},v_{2,m},v^-_{1,m},v^-_{2,m}]$ whose exterior unit normals are
$$
u_{m}^+=\left(\frac{-a_m}{\sqrt{a_m^2+m^2}},0,\frac{m}{\sqrt{a_m^2+m^2}}\right)
\mbox{ \ and \ }
u_{m}^-=\left(\frac{-a_m}{\sqrt{a_m^2+m^2}},0,\frac{-m}{\sqrt{a_m^2+m^2}}\right),
$$
which satisfy $\lim_{m\to\infty}u_{m}^+=u^+$ and $\lim_{m\to\infty}u_{m}^-=u^-$.

For $i=1,2$, we have $h_{\widetilde{P}_m}(u_0)=m$, 
$h_{\widetilde{P}_m}(u_1)=h_{\widetilde{P}_m}(u_2)=\frac{m}{\sqrt{2}}$ and 
 $h_{\widetilde{P}_m}(u_m^+)=h_{\widetilde{P}_m}(u_m^-)=0$, therefore
\begin{eqnarray*}
S_{\widetilde{P}_m,p}(\{u_0\})&=&h_{\widetilde{P}_m}(u_0)^{1-p}\mathcal{H}^2(F_{0,m})=
m^{1-p}8ma_m=8\\
S_{\widetilde{P}_m,p}(\{u_i\})&=&h_{\widetilde{P}_m}(u_i)^{1-p}\mathcal{H}^2(F_{i,m})=
\left(\frac{m}{\sqrt{2}}\right)^{1-p}\sqrt{2}ma_m=2^{\frac{p}2}\mbox{ \ for }i=1,2.
\end{eqnarray*}
Now we translate $\widetilde{P}_m$ in order to alter $S_{\widetilde{P}_m,p}(\{u^+_m\})$. We define $t_m>0$ in a way such that $P_m=\widetilde{P}_m-t_mu_0$ satisfies
$$
h_{P_m}(u^+_m)=m^{\frac{-2}{1-p}}.
$$
It follows that 
$$
m^{\frac{-2}{1-p}}=h_{P_m}(u^+_m)=t_m\langle u^+_m,u_0\rangle=\frac{t_ma_m}{\sqrt{m^2+a_m^2}}
>\frac{t_m}{2m^{3-p}}.
$$
We observe that $r=3-p-\frac2{1-p}<3-2=1$ if $p\in(0,1)$, and hence
$\lim_{m\to\infty}t_m/m=0$. We deduce that
\begin{eqnarray*}
\lim_{m\to\infty}S_{P_m,p}(\{u_0\})&=&8\\
\lim_{m\to\infty}S_{P_m,p}(\{u_i\})&=&2^{\frac{p}2}\mbox{ \ for }i=1,2,\\
\lim_{m\to\infty}S_{P_m,p}(\{u^+_m\})&=&
\lim_{m\to\infty}h_{P_m}(u^+_m)^{1-p}\mathcal{H}^2(F_{m}^+)=
\lim_{m\to\infty}m^{-2}3m\sqrt{m^2+a_m^2}=3.
\end{eqnarray*}
Therefore $S_{P_m,p}$ tends weakly to $\mu$. \ Q.E.D.

\section{Proof of Theorem~\ref{theoplanar0<1} if the measure is concentrated on a closed semi-circle} 

First we show that the $L_p$ surface area measure of a convex body $K$ containing the origin can't be supported on two antipodal points.

\begin{lemma}
\label{SKpcond}
If $K\in{\mathcal K}_{0}^2$, then ${\rm supp}\,S_{K,p}$ is not a pair of antipodal points.
\end{lemma}
\proof We suppose that ${\rm supp}\,S_{K,p}=\{v,-v\}$ for some $v\in S^1$, and seek a contradiction. Let $w\in S^1$ be orthogonal to $v$.

If $o\in{\rm int}\,K$ then ${\rm supp}\,S_{K,p}={\rm supp}\,S_K$, which is not contained in any closed semi-circle. Therefore $o\in\partial K$, and let $C$ be the exterior normal cone at $o$; namely, 
$C\cap S^1=\{u\in S^1:h_K(u)=0\}$. 
Since  ${\rm supp}\,S_{K,p}=\{v,-v\}$, we have $h_K(v)>0$ and $h_K(-v)>0$, and hence 
 $v,-v\not\in C$. Thus  we may assume possibly after replacing $w$ with $-w$ that
$C\cap S^1\subset \Omega(-w,0)$. It follows that $h_K(u)>0$ for $u\in \Omega(w,0)$, and since 
$S_K(\Omega(w,0))>0$,  it also follows that
$$
S_{K,p}(\Omega(w,0))=\int_{\Omega(w,0)}h_K^{1-p}\,dS_K>0.
$$
This contradicts ${\rm supp}\,S_{K,p}=\{v,-v\}$, and proves Lemma~\ref{SKpcond}. \ Q.E.D.\\

Let $\mu$ be a non-trivial measure on $S^1$ that is concentrated on a closed semi-circle 
$\sigma$ of $S^1$
connecting $v,-v\in S^1$ such that ${\rm supp}\,\mu$ is not a pair of antipodal points. 
We may assume that for the $w\in \sigma$ orthogonal to $v$, we have either ${\rm supp}\,\mu=\{w\}$, or
\begin{equation}
\label{winside}
w\in{\rm int}\,{\rm pos}({\rm supp}\,\mu).
\end{equation} 

\noindent {\bf Case 1 } ${\rm supp}\,\mu=\{w\}$

Let $w_1,w_2\in S^1$ such that $w_1+w_2=-w$, and let $K_0$ be the regular triangle 
$$
K_0=\{x\in\R^2:\langle x,w_1\rangle\leq 0,\;\langle x,w_2\rangle\leq 0,\;\langle x,w\rangle\leq 1\}.
$$
For $\lambda=\mu(\{w\})/ S_{K_0,p}(\{w\})$ and $\lambda_0=\lambda^{\frac1{2-p}}$, we have
$S_{\lambda_0K_0,p}=\mu$.\\

\noindent {\bf Case 2 } $w\in{\rm int}\,{\rm pos}({\rm supp}\,\mu)$

Let $A$ be the reflection through the line ${\rm lin}\,v$. We define a measure $\tilde{\mu}$ on $S^1$ by
$$
\tilde{\mu}(\omega)=\mu(\omega)+\mu(A\omega) \mbox{ \ for Borel sets $\omega\subset S^1$}.
$$
We observe that $\tilde{\mu}$ is invariant under $A$,
\begin{eqnarray*}
\tilde{\mu}(\omega)&=&\mu(\omega)\mbox{ \ if $\omega\subset \Omega(w,0)$},\\
\tilde{\mu}(\{v\})&=&2\mu(\{v\}),\\
\tilde{\mu}(\{-v\})&=&2\mu(\{-v\}).
\end{eqnarray*}
 It follows from $w\in{\rm int}\,{\rm pos}({\rm supp}\,\mu)$ that no closed semi-circle contains ${\rm supp}\,\tilde{\mu}$. 
Since the case of Theorem~\ref{theoplanar0<1} when the measure of any open semicircle is positive has been already proved in Section~\ref{secopensemipos}, 
 there exists
a  convex body $\widetilde{K}\in{\mathcal K}_{0}^2$ invariant under $A$ such that $S_{\widetilde{K},p}=\tilde{\mu}$.

We claim that 
\begin{equation}
\label{SKpmu}
S_{K,p}=\mu\mbox{ \ for $K=\{x\in \widetilde{K}:\,\langle x,w\rangle\geq 0\}$}.
\end{equation} 
For any convex body $M$ and $u\in S^1$, we write $F(M,u)=\{x\in M:\,\langle x,u\rangle=h_M(u)\}$ for the face of $M$ with exterior unit normal $u$, and for any $x,y\in \R^2$, we write $[x,y]$ for the convex hull of $x$ and $y$, which is a segment if $x\neq y$. Since $\widetilde{K}$ is invariant under $A$, there exist $t,s\geq 0$ such that
$tv,-sv\in \partial \widetilde{K}$, and the exterior normals at $tv$ and $-sv$ are $v$ and $-v$, respectively.
In addition, $\mathcal{H}^1(F(\widetilde{K},v))=2\,\mathcal{H}^1(F(K,v))$,
$\mathcal{H}^1(F(\widetilde{K},-v))=2\,\mathcal{H}^1(F(K,-v))$ and
$F(K,-w)=[tv,-sv]$.

To prove (\ref{SKpmu}), first we observe that by definition, we have
$$
\mu(\{v\})=\frac{\tilde{\mu}(\{v\})}2=\frac{h_{\widetilde{K}}(v)^{1-p}\cdot \mathcal{H}^1(F(\widetilde{K},v))}2=
h_{K}(v)^{1-p}\cdot \mathcal{H}^1(F(K,v))=S_{K,p}(\{v\}),
$$
and similarly $\mu(\{-v\})=S_{K,p}(\{-v\})$. Next (\ref{SKpdef}) yields that
$$
S_{K,p}(\Omega(-w,0))=\int_{[tv,-sv]}\langle x,w\rangle^{1-p}\,d\mathcal{H}^1(x)
=0=\mu(\Omega(-w,0)).
$$  
Finally, if $\omega\subset \Omega(w,0)$, then $\nu_{\widetilde{K}}^{-1}(\omega)=\nu_{K}^{-1}(\omega)$,
which  yields
$$
\mu(\omega)=\tilde{\mu}(\omega)=S_{\widetilde{K},p}(\omega)=S_{K,p}(\omega),
$$
and in turn (\ref{SKpmu}).

Therefore all we are left to do is to check the symmetries of $\mu$. Actually the only possible symmetry is the reflection $B$ through ${\rm lin}\,w$. In this case, $\tilde{\mu}$ is also invariant under $B$, and hence we may assume that
$\widetilde{K}$ is also invariant under $B$. We conclude that $K$ is invariant under $B$, completing the proof of 
Theorem~\ref{theoplanar0<1}. \ Q.E.D.

\section{Appendix}

Let $p\in(0,1)$, let $\mu$ be a discrete measure on $S^{n-1}$ such that any open hemi-sphere has positive measure, and
let $G\subset O(n)$ is a subgroup such that 
$\mu(\{Au\})=\mu(\{u\})$ for every $u\in S^{n-1}$ and $A\in G$. We review the proof of Theorem~\ref{theopoly0<1} due to Zhu \cite{Zhu15} to show that for the polytope $P$ with $o\in{\rm int}\,P$ and $S_{P,p}=\mu$,  one may 
even assume that $AP=P$ for every $A\in G$.

We set ${\rm supp}\,\mu=\{u_1,\ldots,u_N\}$ and $\alpha_i=\mu(\{u_i\})>0$ for $i=1,\ldots,N$, and we write 
$$
\mathcal{P}^G(u_1,\ldots,u_N)
$$
for the family of $n$-dimensional polytopes whose exterior unit normals are among $u_1,\ldots,u_N$ and which are $G$ invariant. In particular,
 if $P\in \mathcal{P}^G(u_1,\ldots,u_N)$ and $A\in G$, then
$h_P(Au_i)=h_P(u_i)$ for $i=1,\ldots,N$.

 In order to find a polytope $P_0\in \mathcal{P}^G(u_1,\ldots,u_N)$ with $S_{P_0,p}=\mu$, following 
Zhu \cite{Zhu15}, we consider
$$
\Phi_P(\xi)=\int_{S^{n-1}}h_{P-\xi}^p\,d\mu=\sum_{i=1}^N\alpha_i(h_P(u_i)-\langle \xi,u_i\rangle)^p
$$
for $P\in \mathcal{P}^G(u_1,\ldots,u_N)$ and $\xi\in P$, and show that the extremal problem
$$
\inf\left\{\sup_{\xi\in P}\Phi_P(\xi):\,P\in \mathcal{P}^G(u_1,\ldots,u_N)\mbox{ and }V(P)=1\right\}
$$
has a solution that is a dilated copy of $P_0$.

According to Lemma~3.1 and  Lemma 3.2 in  \cite{Zhu15},  if $P\in \mathcal{P}^G(u_1,\ldots,u_N)$, then there exists a unique $\xi(P)\in{\rm int}P$ such that
$$
\sup_{\xi\in P}\Phi_P(\xi)=\Phi_P(\xi(P)).
$$
The uniqueness of $\xi(P)$ yields that
$$
A\xi(P)=\xi(P)\mbox{ \ for $A\in G$}.
$$
We deduce from Lemma 3.3 in  \cite{Zhu15} that $\xi(P)$ is a continuous function of $P$. 

Let $\mathcal{P}_N^G(u_1,\ldots,u_N)$ be the family of all $P\in\mathcal{P}^G(u_1,\ldots,u_N)$ with $N$ facets.
Based on Lemma~3.4 and  Lemma 3.5 in  \cite{Zhu15}, slightly modifying the argument for 
Lemma~3.6 in  \cite{Zhu15}, we deduce the existence of $\widetilde{P}\in \mathcal{P}_N^G(u_1,\ldots,u_N)$
with $V(\widetilde{P})=1$ such that 
$$
\Phi_{\widetilde{P}}(\xi(\widetilde{P}))=
\inf\left\{\Phi_P(\xi(P)):\,P\in \mathcal{P}^G(u_1,\ldots,u_N)\mbox{ and }V(P)=1\right\}.
$$
The only change in the argument for 
Lemma~3.6 in  \cite{Zhu15} is making the definition of $P_\delta$ $G$ invariant. So supposing that
${\rm dim}\,F(\widetilde{P},u_{i_0})\leq n-2$, let $I\subset \{1,\ldots,N\}$ be defined by
$$
\{Au_{i_0}:\,A\in G\}=\{u_i:\,i\in I\}.
$$
Therefore for small $\delta>0$, we set
$$
P_\delta=\{x\in P:\,\langle x,u_i\rangle\leq h_{\widetilde{P}}(u_i)-\delta\mbox{ \ for $i\in I$}\}.
$$
The rest of the argument for Lemma~3.6 in  \cite{Zhu15} carries over.

Finally, in the proof of Theorem~4.1 in  \cite{Zhu15}, the only necessary change is that for the
$\delta_1,\ldots,\delta_N\in \R$ we assume that for any $A\in G$ and $i\in\{1,\ldots,N\}$, if
$u_j=Au_i$, then $\delta_j=\delta_i$.\\

\noindent{\bf Acknowledgement } We are grateful for the referees' helpful comments.

\end{document}